\documentclass[11pt]{book}

\usepackage{amssymb}
\usepackage{amsmath}
\usepackage{amsfonts}
\usepackage{epsfig}

\makeatletter

\renewenvironment{thebibliography}[1]
     {\section*{\bibname}%
      \@mkboth{\MakeUppercase\bibname}{\MakeUppercase\bibname}%
      \list{\@biblabel{\@arabic\c@enumiv}}%
           {\settowidth\labelwidth{\@biblabel{#1}}%
            \leftmargin\labelwidth
            \advance\leftmargin\labelsep
            \@openbib@code
            \usecounter{enumiv}%
            \let\p@enumiv\@empty
            \renewcommand\theenumiv{\@arabic\c@enumiv}}%
      \sloppy
      \clubpenalty4000
      \@clubpenalty \clubpenalty
      \widowpenalty4000%
      \sfcode`\.\@m}
     {\def\@noitemerr
       {\@latex@warning{Empty `thebibliography' environment}}%
      \endlist}

\makeatother

\newcommand{\sect}[1]{\setcounter{equation}{0}\section{#1}}

\renewcommand\bibname{References}

\newtheorem{theorem}{Theorem}[section]
\newtheorem{proposition}[theorem]{Proposition}

\newtheorem{definition}[theorem]{Definition}

\newtheorem{notation}[theorem]{Notation}

\def\cA{{\mathcal A}}

\def\cB{{\mathcal B}}

\def\cF{\mathcal{F}}
\def\cR{\mathcal{R}}

\def\ff{\varphi}

\newcommand{\NN}{ {\mathbb N} }
\newcommand{\RR}{ {\mathbb R} }
\newcommand{\CC}{{\mathbb C}}
\newcommand{\FF}{{\mathbb F}}

\newcommand{\HH}{{\mathcal H}}

\def\kk{\kappa}

\def\tr{\text{\rm tr}}

\newcommand{\distr}{\text{distr}}

\begin{document}
\setlength{\baselineskip}{5.0mm}


\setcounter{chapter}{21}
\chapter[]{Free Probability Theory}
\thispagestyle{empty}

\ \\

\noindent {{\sc Roland Speicher}
\\~\\
Department of Mathematics and Statistics\\
Queen's University\\
Kingston, Ontario K7L 3N6\\
Canada}

\begin{center}
{\bf Abstract}
\end{center}
Free probability theory was created by Dan Voiculescu  around 1985, motivated by his
efforts to understand special classes of von Neumann algebras. His discovery in 1991 that
also random matrices satisfy asymptotically the \emph{freeness} relation transformed the
theory dramatically. Not only did this yield spectacular results about the structure of
operator algebras, but it also brought new concepts and tools into the realm of random
matrix theory. In the following we will give, mostly from the random matrix point of
view, a survey on some of the basic ideas and results of free probability theory.

\sect{Introduction}\label{intro}

Free probability theory allows one to deal with asymptotic eigenvalue distributions in
situations involving several matrices. Let us consider two sequences $A_N$ and $B_N$ of
selfadjoint $N\times N$ matrices such that both sequences have an asymptotic eigenvalue
distribution for $N\to\infty$. We are interested in the asymptotic eigenvalue
distribution of the sequence $f(A_N,B_N)$ for some non-trivial selfadjoint function $f$.
In general, this will depend on the relation between the eigenspaces of $A_N$ and of
$B_N$. However, by the concentration of measure phenomenon, we expect that for large $N$
this relation between the eigenspaces concentrates on \emph{typical} or \emph{generic
positions}, and that then the asymptotic eigenvalue distribution of $f(A_N,B_N)$ depends
in a deterministic way only on the asymptotic eigenvalue distribution of $A_N$ and on the
asymptotic eigenvalue distribution of $B_N$. Free probability theory replaces this vague
notion of \emph{generic position} by the mathematical precise concept of \emph{freeness}
and provides general tools for calculating the asymptotic distribution of $f(A_N,B_N)$
out of the asymptotic distribution of $A_N$ and the asymptotic distribution of $B_N$.

\section[Moment Method and Freeness]
{The Moment Method for Several Random Matrices and the Concept of Freeness}

The \emph{empirical eigenvalue distribution} of a selfadjoint $N\times N$ matrix $A$ is
the probability measure on $\RR$ which puts mass $1/N$ on each of the $N$ eigenvalues
$\lambda_i$ of $A$, counted with multiplicity. If $\mu_A$ is determined by its moments
then it can be recovered from the knowledge of all traces of powers of $A$:
$$\tr(A^k)=\frac 1N\bigl(\lambda_1^k+\cdots+\lambda_N^k\bigr)=\int_\RR t^k d\mu_A(t),$$
where by $\tr$ we denote the normalized trace on matrices (so that we have for the
identity matrix $1$ that $\tr(1)=1$). This is the basis of the \emph{moment method} which
tries to understand the asymptotic eigenvalue distribution of a sequence of matrices by
the determination of the asymptotics of traces of powers.

\begin{definition}
We say that a sequence $(A_N)_{N\in\NN}$ of $N\times N$ matrices \emph{has an asymptotic
eigenvalue distribution} if the limit $\lim_{N\to\infty}\tr(A_N^k)$ exists for all
$k\in\NN$.
\end{definition}

Consider now our sequences $A_N$ and $B_N$, each of which is assumed to have an
asymptotic eigenvalue distribution. We want to understand, in the limit $N\to\infty$, the
eigenvalue distribution of $f(A_N,B_N)$, not just for one $f$, but for a wide class of
different functions. By the moment method, this asks for the investigation of the limit
$N\to\infty$ of $\tr\bigl(f(A_N,B_N)^k\bigr)$ for all $k\in\NN$ and all $f$ in our
considered class of functions. If we choose for the latter all polynomials in
non-commutative variables, then it is clear that the basic objects which we have to
understand in this approach are the asymptotic \emph{mixed moments}
\begin{equation}\label{eq:mixed-moments}
\lim_{N\to\infty}\tr(A_N^{n_1}B_N^{m_1}\cdots A_N^{n_k}B_N^{m_k})\qquad (k\in\NN;\,
n_1,\dots,n_k,m_1,\dots,m_k\in\NN).
\end{equation}

Thus our fundamental problem is the following. If $A_N$ and $B_N$ each have an asymptotic
eigenvalue distribution, and if $A_N$ and $B_N$ are in generic position, do the
asymptotic mixed moments $\lim_{N\to\infty} \tr(A_N^{n_1}B_N^{m_1}\cdots
A_N^{n_k}B_N^{m_k})$ exist? If so, can we express them in a deterministic way in terms of
\begin{equation}\label{eq:A-B-moments}
\left(\lim_{N\to\infty}\tr(\ff(A_N^k))\right)_{k\in\NN}\qquad\text{and}\qquad
\left(\lim_{N\to\infty}\tr(B_N^k)\right)_{k\in\NN}.
\end{equation}

Let us start by looking on the second part of the problem, namely by trying to find a
possible relation between the mixed moments \eqref{eq:mixed-moments} and the moments
\eqref{eq:A-B-moments}. For this we need a simple example of matrix sequences $A_N$ and
$B_N$ which we expect to be in generic position.

Whereas up to now we have only talked about sequences of matrices, we will now go over to
random matrices. Namely, it is actually not clear how to produce two sequences of deterministic
matrices whose eigenspaces are in generic position. However, it is much easier to produce two
such sequences of random matrices for which we have almost surely a generic situation. Indeed,
consider two independent random matrix ensembles $A_N$ and $B_N$, each with almost surely a
limiting eigenvalue distribution, and assume that one of them, say $B_N$, is a \emph{unitarily
invariant ensemble}, which means that the joint distribution of its entries does not change
under unitary conjugation. This implies that taking $U_NB_NU_N^*$, for any unitary $N\times
N$-matrix $U_N$, instead of $B_N$ does not change anything. But then we can use this $U_N$ to
rotate the eigenspaces of $B_N$ against those of $A_N$ into a generic position, thus for
typical realizations of $A_N$ and $B_N$ the eigenspaces should be in a generic position.

The simplest example of two such random matrix ensembles are two independent Gaussian
random matrices $A_N$ and $B_N$. In this case one can calculate everything concretely: in
the limit $N\to\infty$, $\tr(A_N^{n_1}B_N^{m_1}\cdots A_N^{n_k}B_N^{m_k})$ is almost
surely given by the number of non-crossing or planar pairings of the pattern
$$\underbrace{A\cdot A\cdots A}_{\text{$n_1$-times}}\cdot
\underbrace{B\cdot B\cdots B}_{\text{$m_1$-times}}\cdots \underbrace{A\cdot A\cdots
A}_{\text{$n_k$-times}}\cdot \underbrace{B\cdot B\cdots B}_{\text{$m_k$-times}},
$$
which do not pair $A$ with $B$. (A pairing is a decomposition of the pattern into pairs
of letters; if we connect the two elements from each pair by a line, drawn in the
half-plane below the pattern, then non-crossing means that we can do this without getting
crossings between lines for different pairs.)

After some contemplation, it becomes obvious that this implies that the trace of a
corresponding product of centered powers,
\begin{multline}\label{eq:freeness-rm}
\lim_{N\to\infty}\tr\Big(\bigl(A_N^{n_1}-\lim_{M\to\infty}\tr(A_M^{n_1})\cdot1\bigr)\cdot
\bigl(B_N^{m_1}-\lim_{M\to\infty}\tr(B_M^{m_1})\cdot1\bigr)\cdots
\\ \cdots \bigl(A_M^{n_k}-\lim_{M\to\infty}\tr(A_M^{n_k})\cdot1\bigr)
\cdot\bigl(B_N^{m_k}-\lim_{M\to\infty}\tr(B_M^{m_k})\cdot1\bigr)\Big)
\end{multline}
is given by the number of non-crossing pairings which do not pair $A$ with $B$ and for
which, in addition, each group of $A$'s and each group of $B$'s is connected with some
other group. It is clear that if we want to connect the groups in this way we will get
some crossing between the pairs, thus there are actually no pairings of the required form
and we have that the term \eqref{eq:freeness-rm} is equal to zero.

One might wonder what advantage is gained by trading the explicit formula for mixed
moments of independent Gaussian random matrices for the implicit relation
\eqref{eq:freeness-rm}? The drawback to the explicit formula for mixed moments of
independent Gaussian random matrices is that the asymptotic formula for
$\tr(A_N^{n_1}B_N^{m_1}\cdots A_N^{n_k}B_N^{m_k})$ will be different for different random
matrix ensembles (and in many cases an explicit formula fails to exist). However, the
vanishing of \eqref{eq:freeness-rm} remains valid for many matrix ensembles. The
vanishing of \eqref{eq:freeness-rm} gives a precise meaning to our idea that the random
matrices should be in generic position; it constitutes Voiculescu's definition of
asymptotic freeness.

\begin{definition} Two sequences of matrices $(A_N)_{N\in\NN}$ and
$(B_N)_{N\in\NN}$ are \emph{asymptotically free} if we have the vanishing of
\eqref{eq:freeness-rm} for all $k\geq 1$ and all $n_1,m_1,$ $\dots$,$ n_k,m_k\geq 1$.
\end{definition}

Provided with this definition, the intuition that unitarily invariant random matrices should
give rise to generic situations becomes now a rigorous theorem. This basic observation was
proved by Voiculescu \cite{Voi91} in 1991.

\begin{theorem}\label{thm:random-matrix-freeness}
Consider $N\times N$ random matrices $A_N$ and $B_N$ such
that: both $A_N$ and $B_N$ have almost surely an asymptotic eigenvalue distribution for
$N\to\infty$; $A_N$ and $B_N$ are independent; $B_N$ is a unitarily invariant ensemble.
Then, $A_N$ and $B_N$ are almost surely asymptotically free.
\end{theorem}

In order to prove this, one can replace $A_N$ and $B_N$ by $A_N$ and $U_NB_NU_N^*$, where
$U_N$ is a Haar unitary random matrix (i.e., from the ensemble of unitary matrices
equipped with the normalized Haar measure as probability measure); furthermore, one can
then restrict to the case where $A_N$ and $B_N$ are deterministic matrices. In this form
it reduces to showing almost sure asymptotic freeness between Haar unitary matrices and
deterministic matrices. The proof of that statement proceeds then as follows. First one
shows asymptotic freeness in the mean and then one strengthens this to almost sure
convergence.

The original proof of Voiculescu \cite{Voi91} for the first step reduced the asymptotic
freeness for Haar unitary matrices to a corresponding statement for non-selfadjoint
Gaussian random matrices; by realizing the Haar measure on the group of unitary matrices
as the pushforward of the Gaussian measure under taking the phase. The asymptotic
freeness result for Gaussian random matrices can be derived quite directly by using the
genus expansion for their traces. Another more direct way to prove the averaged version
of unitary freeness for Haar unitary matrices is due to Xu \cite{Xu97} and relies on
\emph{Weingarten type formulas} for integrals over products of entries of Haar unitary
matrices.

In the second step, in order to strengthen the above result to almost sure asymptotic
freeness one can either \cite{Voi91} invoke concentration of measure results of Gromov
and Milman (applied to the unitary group) or \cite{Spe93, Hia00} more specific estimates
for the variances of the considered sequence of random variables.

Though unitary invariance is the most intuitive reason for having asymptotic freeness among
random matrices, it is not a necessary condition. For example, the above theorem includes the
case where $B_N$ are Gaussian random matrices. If we generalize those to Wigner matrices (where
the entries above the diagonal are i.i.d, but not necessarily Gaussian), then we loose the
unitary invariance, but the conclusion of the above theorem still holds true. More precisely,
we have the following theorem.

\begin{theorem}
Let $X_N$ be a selfadjoint Wigner matrix, such that the distribution of the entries is
centered and has all moments, and let $A_N$ be a random matrix which is independent from
$X_N$. If $A_N$ has almost surely an asymptotic eigenvalue distribution and if we have
$$\sup_{N\in \NN} \Vert A_N\Vert<\infty,$$
then $A_N$ and $X_N$ are almost surely asymptotically free.
\end{theorem}

The case where $A_N$ consists of block diagonal matrices was treated by Dykema \cite{Dyk93},
for the general version see \cite{MS, And10}.

\sect{Basic Definitions}

The freeness relation, which holds for many random matrices asymptotically, was actually
discovered by Voiculescu in a quite different context; namely canonical generators in
operator algebras given in terms of free groups satisfy the same relation with respect to
a canonical state, see Section \ref{sect:freegroup}. Free probability theory investigates
these freeness relations abstractly, inspired by the philosophy that freeness should be
considered and treated as a kind of non-commutative analogue of the classical notion of
independence.

Some of the main probabilistic notions used in free probability are the following.

\begin{notation}
A pair $(\cA,\ff)$ consisting of a unital algebra $\cA$ and a linear functional
$\ff:\cA\to\CC$ with $\ff(1)=1$ is called a \emph{non-commutative probability space}.
Often the adjective ``non-commutative'' is just dropped. Elements from $\cA$ are
addressed as \emph{(non-commutative) random variables}, the numbers $\ff(a_{i(1)}\cdots
a_{i(n)})$ for such random variables $a_1,\dots,a_k\in\cA$ are called \emph{moments}, the
collection of all moments is called the \emph{joint distribution of $a_1,\dots,a_k$}.
\end{notation}

\begin{definition}\label{def:freeness}
Let $(\cA,\ff)$ be a non-commutative probability space and let $I$ be an index set.

1) Let, for each $i\in I$, $\cA_i\subset \cA$, be a unital subalgebra. The subalgebras
$(\cA_i)_{i\in I}$ are called \emph{free} or \emph{freely independent}, if $\ff(a_1\cdots
a_k)=0$ whenever we have: $k$ is a positive integer; $a_j\in\cA_{i(j)}$ (with $i(j)\in
I$) for all $j=1,\dots,k$; $\ff(a_j)=0$ for all $j=1,\dots,k$; and neighboring elements
are from different subalgebras, i.e., $i(1)\not=i(2), i(2)\not= i(3),\dots,
i(k-1)\not=i(k)$.

2) Let, for each $i\in I$, $a_i\in\cA$. The elements $(a_i)_{i\in I}$ are called
\emph{free} or \emph{freely independent}, if their generated unital subalgebras are free,
i.e., if $(\cA_i)_{i\in I}$ are free, where, for each $i\in I$, $\cA_i$ is the unital
subalgebra of $\cA$ which is generated by $a_i$.
\end{definition}

Freeness, like classical independence, is a rule for calculating mixed moments from knowledge
of the moments of individual variables. Indeed, one can easily show by induction that if
$(\cA_i)_{i\in I}$ are free with respect to $\ff$, then $\ff$ restricted to the algebra
generated by all $\cA_i$, $i\in I$, is uniquely determined by $\ff\vert_{\cA_i}$ for all $i\in
I$ and by the freeness condition. For example, if $\cA$ and $\cB$ are free, then one has for
$a,a_1,a_2\in\cA$ and $b,b_1,b_2\in\cB$ that $\ff(ab)=\ff(a)\ff(b)$,
$\ff(a_1ba_2)=\ff(a_1a_2)\ff(b)$, and
$\ff(a_1b_1a_2b_2)=\ff(a_1a_2)\ff(b_1)\ff(b_2)+\ff(a_1)\ff(a_2)\ff(b_1b_2)-
\ff(a_1)\ff(b_1)\ff(a_2) \ff(b_2)$. Whereas the first two factorizations are the same as for
the expectation of independent random variables, the last one is different, and more
complicated, from the classical situation. It is important to note that freeness plays a
similar role in the non-commutative world as independence plays in the classical world, but
that freeness is not a generalization of independence: independent random variables can be free
only in very trivial situations. Freeness is a theory for non-commuting random variables.

\sect{Combinatorial Theory of Freeness}

The defining relations for freeness from Def. \ref{def:freeness} are quite implicit and not
easy to handle directly. It has turned out that replacing moments by other quantities,
so-called \emph{free cumulants}, is advantageous for many questions. In particular, freeness is
much easier to describe on the level of free cumulants. The relation between moments and
cumulants is given by summing over non-crossing partitions. This combinatorial theory of
freeness is due to Speicher \cite{Spe94}; many consequences of this approach were worked out by
Nica and Speicher, see \cite{NS}.

\begin{definition}
For a unital linear functional $\ff:\cA\to\CC$ on a unital algebra $\cA$ we define
\emph{cumulant functionals} $\kk_n:\cA^n\to\CC$ (for all $n\geq 1$) by the
\emph{moment-cumulant relations}
\begin{equation}\label{eq:moment-cumulant}
\ff(a_1\cdots a_n)=\sum_{\pi\in NC(n)}\kk_\pi[a_1,\dots, a_n].
\end{equation}
\end{definition}

In equation \eqref{eq:moment-cumulant} the summation is running over \emph{non-crossing
partitions} of the set $\{a_1,a_2,\dots,a_n\}$; those are decompositions of that set into
disjoint non-empty subsets, called \emph{blocks}, such that there are no crossings
between different blocks. In diagrammatic terms this means that if we draw the blocks of
such a $\pi$ below the points $a_1,a_2,\dots, a_n$, then we can do this without having
crossings in our picture. The contribution $\kk_\pi$ in \eqref{eq:moment-cumulant} of
such a non-crossing $\pi$ is a product of cumulants corresponding to the block structure
of $\pi$. For each block of $\pi$ we have as a factor a cumulant which contains as
arguments those $a_i$ which are connected by that block.

An example of a non-crossing partition $\pi$ for $n=10$ is
\setlength{\unitlength}{0.5cm}
$$\begin{picture}(9,4)\thicklines \put(0,0){\line(0,1){3}} \put(0,0){\line(1,0){9}}
\put(9,0){\line(0,1){3}} \put(1,1){\line(0,1){2}} \put(1,1){\line(1,0){7}}
\put(4,1){\line(0,1){2}} \put(8,1){\line(0,1){2}} \put(2,2){\line(0,1){1}}
\put(2,2){\line(1,0){1}} \put(3,2){\line(0,1){1}} \put(5,2){\line(0,1){1}}
\put(6,2){\line(0,1){1}} \put(6,2){\line(1,0){1}} \put(7,2){\line(0,1){1}}
\put(-0.1,3.3){$a_1$} \put(0.9,3.3){$a_2$} \put(1.9,3.3){$a_3$} \put(2.9,3.3){$a_4$}
\put(3.9,3.3){$a_5$} \put(4.9,3.3){$a_6$} \put(5.9,3.3){$a_7$} \put(6.9,3.3){$a_8$}
\put(7.9,3.3){$a_9$} \put(8.7,3.3){$a_{10}$}
\end{picture}$$
In this case the blocks are $\{a_1,a_{10}\}$, $\{a_2,a_5,a_9\}$, $\{a_3,a_4\}$, $\{a_6\}$, and
$\{a_7,a_8\}$; and the corresponding contribution $\kk_\pi$ in \eqref{eq:moment-cumulant} is
given by
$$
\kk_\pi[a_1,\dots,a_{10}]=\kk_2(a_1,a_{10}) \cdot \kk_3\bigl(a_2,a_5,a_9) \cdot \kk_2(a_3,a_4)
\cdot\kk_1(a_6)\cdot \kk_2(a_7,a_8).
$$

Note that in general there is only one term in \eqref{eq:moment-cumulant} involving the highest
cumulant $\kk_n$, thus the moment cumulant formulas can inductively be resolved for the $\kk_n$
in terms of the moments. More concretely, the set of non-crossing partitions forms a lattice
with respect to refinement order and the $\kk_n$ are given by the \emph{M\"obius inversion} of
the formula \eqref{eq:moment-cumulant} with respect to this order.

For $n=1$, we get the mean, $\kk_1(a_1)=\ff(a_1)$ and for $n=2$ we have the covariance,
$\kk_2(a_1,a_2)=\ff(a_1a_2)-\ff(a_1)\ff(a_2)$.

The relevance of the $\kk_n$ in our context is given by the following characterization of
freeness.

\begin{theorem}
Freeness is equivalent to the vanishing of mixed cumulants. More precisely, the fact that
$(a_i)_{i\in I}$ are free is equivalent to: $\kk_n(a_{i(1)},\dots,a_{i(n)})=0$ whenever
$n\geq 2$ and there are $k,l$ such that $i(k)\not=i(l)$.
\end{theorem}

This description of freeness in terms of free cumulants is related to the planar
approximations in random matrix theory. In a sense some aspects of this theory of
freeness were anticipated (but mostly neglected) in the physics community in the paper
\cite{Cvi82}.

\sect{Free Harmonic Analysis}

For a meaningful harmonic analysis one needs some positivity structure for the
non-commutative probability space $(\cA,\ff)$. We will usually consider selfadjoint
random variables and $\ff$ should be positive. Formally, a good frame for this is a
\emph{$C^*$-probability space}, where $\cA$ is a $C^*$-algebra (i.e., a norm-closed
$*$-subalgebra of the algebra of bounded operators on a Hilbert space) and $\ff$ is a
state, i.e. it is positive in the sense $\ff(aa^*)\geq 0$ for all $a\in \cA$. Concretely
this means that our random variables can be realized as bounded operators on a Hilbert
space and $\ff$ can be written as a vector state $\ff(a)=\langle a\xi,\xi\rangle$ for
some unit vector $\xi$ in the Hilbert space.

In such a situation the distribution of a selfadjoint random variable $a$ can be identified
with a compactly supported probability measure $\mu_a$ on $\RR$, via
$$\ff(a^n)=\int_\RR t^nd\mu_a(t)\qquad\text{for all $n\in\NN$}.$$

\subsection{Sums of free variables: the $\cR$-transform}

Consider two selfadjoint random variables $a$ and $b$ which are free. Then, by freeness, the
moments of $a+b$ are uniquely determined by the moments of $a$ and the moments of $b$.

\begin{notation}
We say the distribution of $a+b$ is the \emph{free convolution}, denoted by $\boxplus$, of the
distribution of $a$ and the distribution of $b$,
$$\mu_{a+b}=\mu_a\boxplus \mu_b.$$
\end{notation}

\begin{notation}\label{not:Cauchy}
For a random variable $a$ we define its \emph{Cauchy transform $G$} and its
\emph{$\cR$-transform $\cR$} by
$$G(z)=\frac 1z+\sum_{n=1}^\infty \frac{\ff(a^n)} {z^{n+1}}\qquad\text{and}\qquad
\cR(z)=\sum_{n=1}^\infty \kk_n(a,\dots,a) z^{n-1}.
$$
\end{notation}
One can see quite easily that the moment-cumulant relations \eqref{eq:moment-cumulant} are
equivalent to the following functional relation
\begin{equation}\label{eq:R-G}
\frac 1{G(z)}+\cR(G(z))=z.
\end{equation}
Combined with the additivity of free cumulants under free convolution, which follows
easily by the vanishing of mixed cumulants in free variables, this yields the following
basic theorem of Voiculescu.

\begin{theorem}\label{thm:R-transform}
Let $G(z)$
be the Cauchy-transform of $a$, as defined in Notation \ref{not:Cauchy} and define its
$\cR$-transform by the relation \eqref{eq:R-G}. Then we have
$$\cR^{a+b}(z)=\cR^a(z)+\cR^b(z)$$
if $a$ and $b$ are free.
\end{theorem}
We have defined the Cauchy and the $\cR$-transform here only as formal power series. Also
\eqref{eq:R-G} is proved first as a relation between formal power series. But if $a$ is a
selfadjoint element in a $C^*$-probability space, then $G$ is also the analytic function
$$G:\CC^+ \to \CC^-; \qquad G(z)=\ff\left(\frac 1{z-a}\right)=\int_\RR\frac 1{z-t}d\mu_a(t);$$
and one can also show that \eqref{eq:R-G} defines then $\cR$ as an analytic function on a
suitably chosen subset of $\CC^+$. In this form Theorem \ref{thm:R-transform} is amenable to
analytic manipulations and so gives an effective algorithm for calculating free convolutions.
This can be used to calculate the asymptotic eigenvalue distribution of sums of random matrices
which are asymptotically free.

Furthermore, by using analytic tools around the Cauchy transform (which exist for any
probability measure on $\RR$) one can extend the definition of and most results on free
convolution to all probability measures on $\RR$. See \cite{Ber93,Voi00} for more
details.

We would like to remark that the machinery of free convolution was also found around the
same time, independently from Voiculescu and independently from each other, by different
researchers in the context of random walks on the free product of groups: by Woess, by
Cartwright and Soardi, and by McLaughlin; see, for example, \cite{Woe86}.

\subsection{Products of free variables: the $S$-transform}

Consider $a$, $b$ free. Then, by freeness, the moments of $ab$ are uniquely determined by the
moments of $a$ and the moments of $b$.

\begin{notation}
We say the distribution of $ab$ is the \emph{free multiplicative convolution}, denoted by
$\boxtimes$, of the distribution of $a$ and the distribution of $b$,
$$\mu_{ab}=\mu_a\boxtimes \mu_b.$$
\end{notation}

Note: even if we start from selfadjoint $a$ and $b$, their product $ab$ is not selfadjoint,
unless $a$ and $b$ commute (which is rarely the case, when $a$ and $b$ are free). Thus the
above does not define an operation on probability measures on $\RR$ in general. However, if one
of the operators, say $a$, is positive (and thus $\mu_a$ supported on $\RR_+$), then $a^{1/2} b
a^{1/2}$ makes sense; since it has the same moments as $ab$ (note for this that the relevant
state is a trace, as the free product of traces is tracial) we can identify $\mu_{ab}$ then
with the probability measure $\mu_{a^{1/2} b a^{1/2}}$.

Again, Voiculescu introduced an analytic object which allows to deal effectively with
this multiplicative free convolution.

\begin{theorem}\label{thm:S-transform}
Put $M_a(z):=\sum_{m=0}^\infty \ff(a^m) z^m$ and define the \emph{{$S$}-transform of {$a$}} by
$$S_a(z):=\frac{1+z}z M_a^{<-1>}(z),$$
where $M^{<-1>}$ denotes the inverse of $M$ under composition. Then we have
$$S_{ab}(z)=S_a(z)\cdot S_b(z)$$
if $a$ and $b$ are free.
\end{theorem}

As in the additive case, the moment generating series $M$ and the $S$-transform are not just
formal power series, but analytic functions on suitably chosen domains in the complex plane.
For more details, see \cite{Ber93,Hia00}.

\subsection{The free central limit theorem}

One of the first theorems in free probability theory, proved by Voiculescu in 1985, was
the free analogue of the central limit theorem. Surprisingly, it turned out that the
analogue of the Gaussian distribution in free probability theory is the semicircular
distribution.

\begin{definition}
Let $(\cA,\ff)$ be a $C^*$-probability space. A selfadjoint element $s\in\cA$ is called
\emph{semicircular} (of variance 1) if its distribution $\mu_s$ is given by the
probability measure with density $\frac 1{2\pi}\sqrt{4-t^2}$ on the interval $[-2,+2]$.
Alternatively, the moments of $s$ are given by the Catalan numbers,
$$\ff(s^n)=\begin{cases}
\frac 1{k+1}\binom{2k}k,& \text{if $n=2k$ even}\\
0,&\text{if $n$ odd}
\end{cases}$$
\end{definition}

\begin{theorem}
If $\nu$ is a compactly supported probability measure on $\RR$ with vanishing mean and
variance 1, then
$$D_{1/\sqrt N}\nu^{\boxplus N} \Rightarrow \mu_s,$$
where $D_\alpha$ denotes the dilation of a measure by the factor $\alpha$, and
$\Rightarrow$ means weak convergence.
\end{theorem}
By using the analytic theory of $\boxplus$ for all, not necessarily compactly supported,
probability measures on $\RR$, the free central limit theorem can also be extended to
this general situation.

The occurrence of the semicircular distribution as limit both in Wigner's semicircle law
as well as in the free central limit theorem was the first hint of a relationship between
free probability theory and random matrices. The subsequent development of this
connection culminated in Voiculescu's discovery of asymptotic freeness between large
random matrices, as exemplified in Theorem \ref{thm:random-matrix-freeness}. When this
contact was made between freeness and random matrices, the previously introduced $\cR$-
and $S$-transforms gave powerful new techniques for calculating asymptotic eigenvalue
distributions of random matrices. For computational aspects of these techniques we refer
to \cite{Rao09}, for applications in electrical engineering see \cite{Tul04}, and also
Chapter 40.

\begin{figure}[h]
\unitlength1cm
\begin{center}
\begin{picture}(15,4)
\centerline{\epsfig{figure=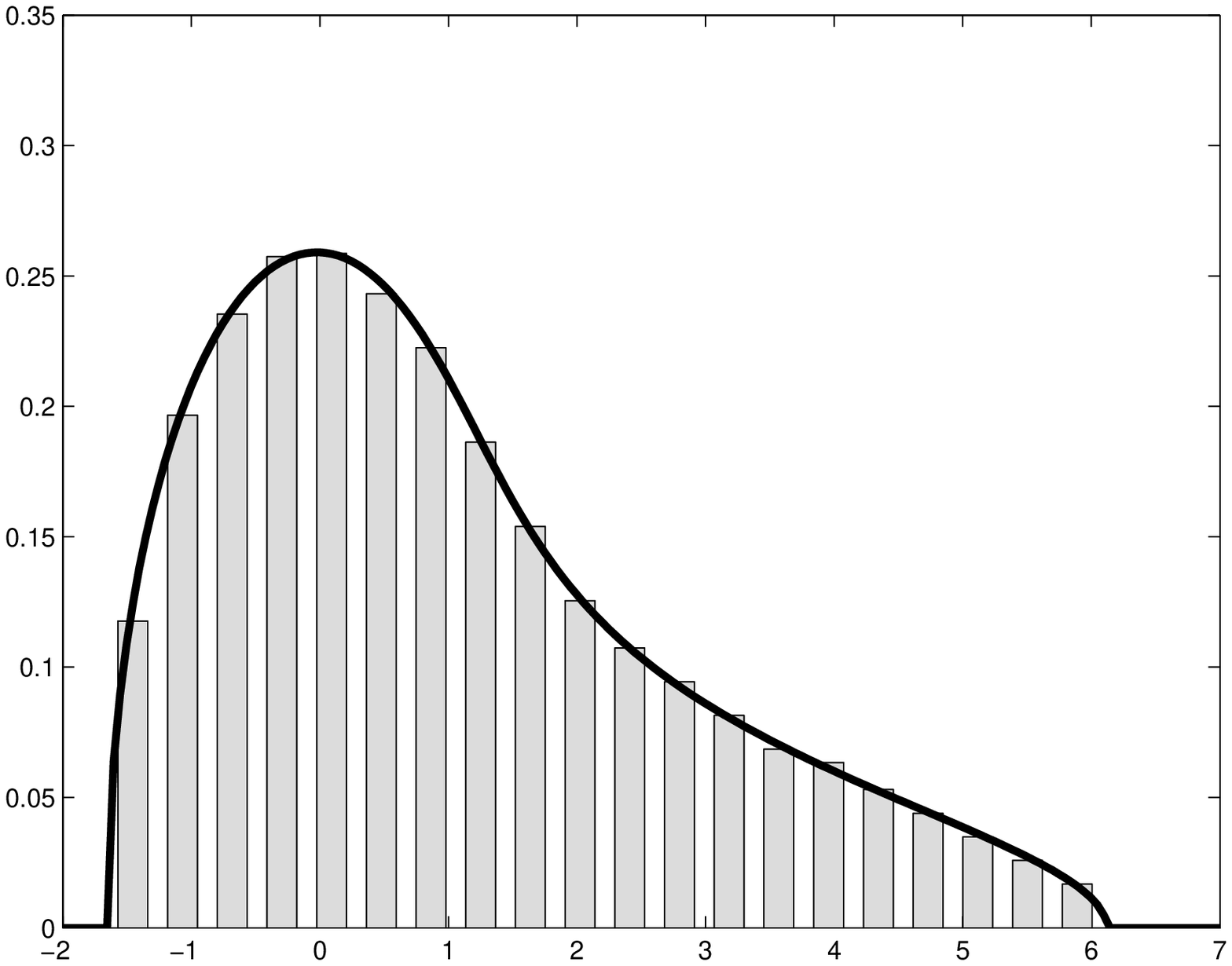,width=11pc},\qquad
\epsfig{figure=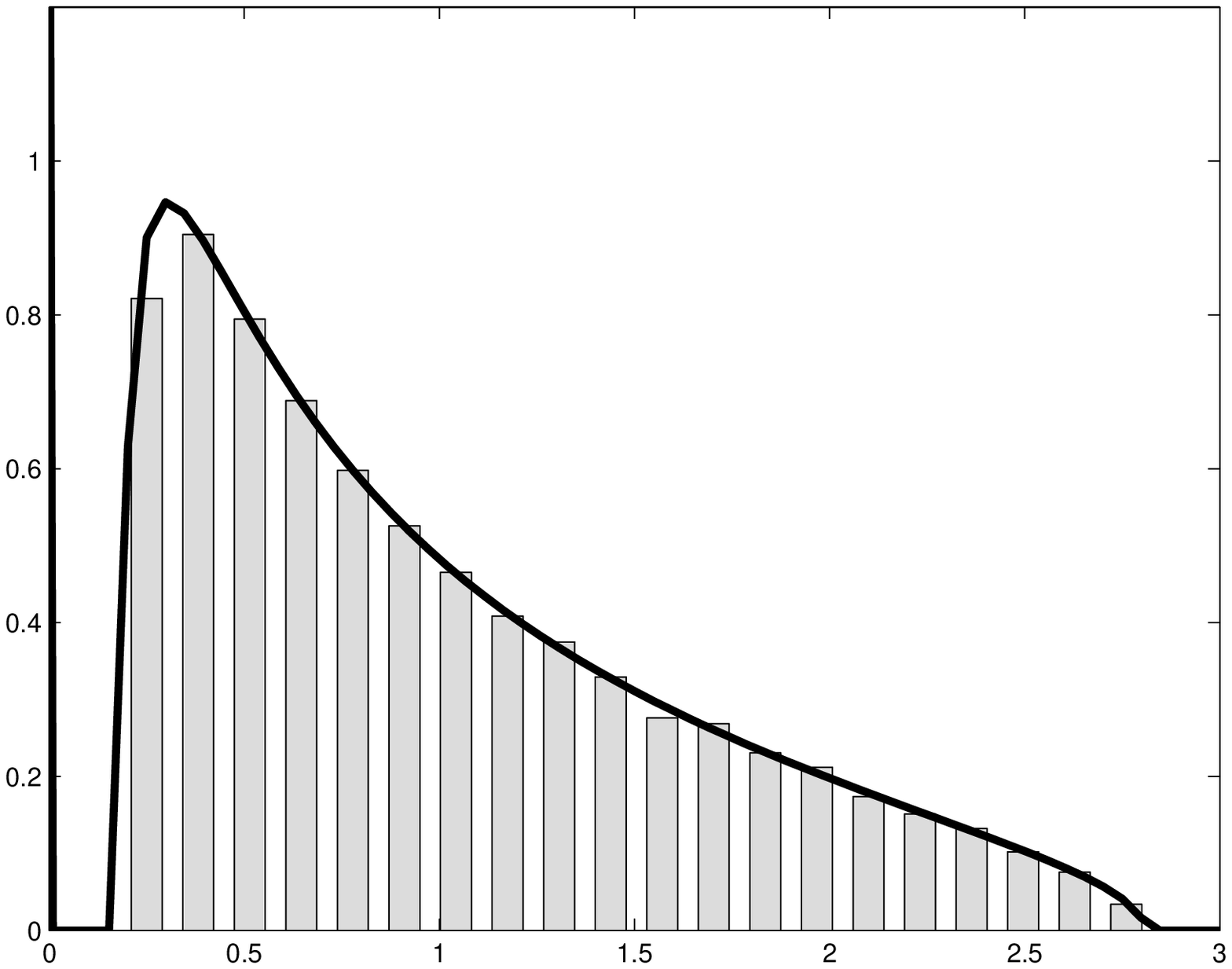,width=11pc},}
\end{picture}
\end{center}
\caption{Comparison of free probability result with histogram of eigenvalues of an $N\times N$
random matrix, for $N=2000$: (i) histogram of the sum of independent Gaussian and Wishart
matrices, compared with the free convolution of semicircular and free Poisson distribution
(rate $\lambda=1/2$), calculated by using the $\cR$-transform; (ii) histogram of the product of
two independent Wishart matrices, compared with the free multiplicative convolution of two free
Poisson distributions (both with rate $\lambda=5$), calculated by using the $S$-transform
\label{fig:26-two}}
\end{figure}

\subsection{Free Poisson distribution and Wishart matrices}

There exists a very rich free parallel of classical probability theory, of which the free
central limit theorem is just the starting point. In particular, one has the free
analogue of infinitely divisible and of stable distributions and corresponding limit
theorems. For more details and references see \cite{Ber99,Voi00}.

Let us here only present as another instance of this theory the free Poisson
distribution. As with the semicircle distribution, the free counterpart of the Poisson
law, which is none other than the Marchenko-Pastur distribution, appears very naturally
as the asymptotic eigenvalue distribution of an important class of random matrices,
namely Wishart matrices.

As in the classical theory the Poisson distribution can be described by a limit theorem. The
following statement deals directly with the more general notion of a compound free Poisson
distribution.

\begin{proposition}\label{prop:limit-compound-Poisson}
Let $\lambda\geq 0$ and $\nu$ a probability measure on $\RR$ with compact support. Then
the weak limit for $N\to\infty$ of
$$\left(
\Bigl(1-\frac \lambda N\Bigr)\delta_0+\frac \lambda N\nu\right)^{ \boxplus N}$$ has free
cumulants $(\kk_n)_{n\geq 1}$ which are given by $\kk_n=\lambda \cdot m_n(\nu)$ ($n\geq
1$) ($m_n$ denotes here the $n$-th moment) and thus an $\cR$-transform of the form
$$\cR(z)=\lambda\int_\RR\frac x{1-xz}\,d\nu(x) .$$
\end{proposition}

\begin{definition}\label{def:compound-Poisson}
The probability measure appearing in the limit of Prop. \ref{prop:limit-compound-Poisson}
is called a \emph{compound free Poisson distribution} with rate $\lambda$ and jump
distribution $\nu$.
\end{definition}

Such compound free Poisson distributions show up in the random matrix context as follows.
Consider: rectangular Gaussian $M\times N$ random matrices $X_{M,N}$, where all entries are
independent and identically distributed according to a normal distribution with mean zero and
variance $1/N$; and a sequence of deterministic $N\times N$ matrices $T_N$ such that the
limiting eigenvalue distribution $\mu_T$ of $T_N$ exists. Then almost surely, for
$M,N\to\infty$ such that $N/M\to\lambda$, the limiting eigenvalue distribution of
$X_{M,N}T_NX_{M,N}^*$ exists, too, and it is given by a compound free Poisson distribution with
rate $\lambda$ and jump distribution $\mu_T$.

One notes that the above frame of rectangular matrices does not fit directly into the theory
presented up to now (so it is, e.g., not clear what asymptotic freeness between $X_{M,N}$ and
$T_N$ should mean). However, rectangular matrices can be treated in free probability by either
embedding them into bigger square matrices and applying some compressions at appropriate stages
or, more directly, by using a generalization of free probability due to Benaych-Georges
\cite{Ben09} which is tailor-made to deal with rectangular random matrices.

\sect{Second Order Freeness}

Asymptotic freeness of random matrices shows that the mixed moments of two ensembles in generic
position are deterministically calculable from the moments of each individual ensemble. The
formula for the calculation of the mixed moments is the essence of the concept of freeness. The
same philosophy applies also to finer questions about random matrices, most notably to
\emph{global fluctuations of linear statistics}. With this we mean the following: for many
examples (like Gaussian or Wishart) of random matrices $A_N$ the magnified fluctuations of
traces around the limiting value, $N\bigl(\tr(A_N^k)-\lim_{M\to\infty}\tr(A_M^k)\bigr)$, form
asymptotically a Gaussian family. If we have two such ensembles in generic position (e.g., if
they are independent and one of them is unitarily invariant), then this is also true for mixed
traces and the covariance of mixed traces is determined in a deterministic way by the
covariances for each of the two ensembles separately. The formula for the calculation of the
mixed covariances constitutes the definition of the concept of \emph{second order freeness}.
There exist again cumulants and an $R$-transform on this level, which allow explicit
calculations. For more details, see \cite{Col07,MS}.

\sect{Operator-Valued Free Probability Theory}

There exists a generalization of free probability theory to an \emph{operator-valued} level,
where the complex numbers $\CC$ and the expectation state $\ff:\cA\to\CC$ are replaced by an
arbitrary algebra $\cB$ and a conditional expectation $E:\cA\to\cB$. The formal structure of
the theory is very much the same as in the scalar-valued case, one only has to take care of the
fact that the ``scalars'' from $\cB$ do not commute with the random variables.

\begin{definition}
1) Let $\cA$ be a unital algebra and consider a unital subalgebra $\cB\subset \cA$. A linear
map $E:\cA\to\cB$ is a \emph{conditional expectation} if $E[b]=b$ for all $b\in\cB$ and
$E[b_1ab_2]=b_1E[a]b_2$ for all $a\in\cA$ and all $b_1,b_2\in\cB$. An \emph{operator-valued
probability space} consists of $\cB\subset \cA$ and a conditional expectation $E:\cA\to\cB$.

2) Consider an operator-valued probability space $\cB\subset\cA$, $E:\cA\to\cB$. Random
variables $(x_i)_{i\in I} \subset\cA$ are \emph{free with respect to $E$} or \emph{free
with amalgamation over $\cB$} if $E[p_1(x_{i(1)})\cdots p_k(x_{i(k)})]=0$, whenever
$k\in\NN$, $p_j$ are elements from the algebra generated by $\cB$ and $x_{i(j)}$,
neighboring elements are different, i.e., $i(1)\not=i(2)\not=\cdots\not= i(k)$,  and we
have $E[p_j(x_{i(j)}]=0$ for all $j=1,\dots,k$.
\end{definition}

Voiculescu introduced this operator-valued version of free probability theory in \cite{Voi85}
and provided in \cite{Voi95} also a corresponding version of free convolution and
$\cR$-transform. A combinatorial treatment was given by Speicher \cite{Spe98} who showed that
the theory of free cumulants has also a nice counterpart in the operator-valued frame.

For $a\in\cA$ we define its \emph{(operator-valued) Cauchy transform} $G_a:\cB\to\cB$ by
$$G_a(b):=E[\frac 1{b-a}]=\sum_{n\geq 0} E[b^{-1}(ab^{-1})^n].$$
The \emph{operator-valued $\cR$-transform} of $a$, $\cR_a:\cB\to\cB$, can be defined as a
power series in operator-valued free cumulants, or equivalently by the relation $b
G(b)=1+\cR(G(b)) \cdot G(b)$ or $G(b)=\bigl({b-\cR(G(b))}\bigr)^{-1}$. One has then as
before: If $x$ and $y$ are free over $\cB$, then $\cR_{x+y}(b)=\cR_x(b)+\cR_y(b)$.
Another form of this is the \emph{subordination property}
$G_{x+y}(b)=G_x\bigl[b-\cR_y\bigl(G_{x+y}(b)\bigr)\bigr]$.

There exists also the notion of a \emph{semicircular element} $s$ in the operator-valued
world. It is characterized by the fact that only its second order free cumulants are
different from zero, or equivalently that its $\cR$-transform is of the form
$\cR_s(b)=\eta(b)$, where $\eta:\cB\to\cB$ is the linear map given by $\eta(b)=E[sbs]$.
Note that in this case the equation for the Cauchy transform reduces to
\begin{equation}\label{eq:G-equation}
b G(b)=1+\eta[G(b)]\cdot G(b);
\end{equation}
more generally, if we add an $x\in\cB$, for which we have
$G_x(b)=E[(b-x)^{-1}]=(b-x)^{-1}$, we have for the Cauchy transform of $x+s$ the implicit
equation
\begin{equation}\label{eq:x+s}
G_{x+s}(b)=G_x\bigl[b-\cR_s\bigl(G_{x+s}(b)\bigr)\bigr]=\bigl(b-\eta[G_{x+s}(b)]
-x\bigr)^{-1}.
\end{equation}

It was observed by Shlyakhtenko \cite{Shl} that operator-valued free probability theory
provides the right frame for dealing with more general kind of random matrices. In
particular, he showed that so-called \emph{band matrices} become asymptotically
operator-valued semicircular elements over the limit of the diagonal matrices.

\begin{theorem}
Suppose that $A_N=A_N^*$ is an $N\times N$ random band matrix, i.e.,
$A_N=(a_{ij})_{i,j=1}^N$, where $\{a_{ij}\mid i\leq j\}$ are centered independent complex
Gaussian random variables, with
$E[a_{ij}\overline{a}_{ij}]=(1+\delta_{ij}\sigma^2(i/N,j/N))/N$ for some $\sigma^2\in
L^\infty([0,1]^2)$. Let $\cB_N$ be the diagonal $N\times N$ matrices, and embed $\cB_N$
into $\cB:=L^\infty[0,1]$ as step functions. Let $B_N\in\cB_N$ be selfadjoint diagonal
matrices such that $B_N\to f\in L^\infty[0,1]$ in $\Vert \cdot \Vert_\infty$. Then the
limit distribution of $B_N+A_N$ exists, and its Cauchy transform $G$ is given by
$$G(z)=\int_0^1 g(z,x)dx,$$
where $g(z,x)$ is analytic in $z$ and satisfies
\begin{equation}\label{eq:gzx}
g(z,x)=\left[z-f(x)-\int_0^1\sigma^2(y,x)g(z,y)dy\right]^{-1}.
\end{equation}
\end{theorem}

Note that \eqref{eq:gzx} is nothing but the general equation \eqref{eq:x+s} specified to
the situation $\cB=L^\infty[0,1]$ and $\eta:L^\infty[0,1]\to L^\infty[0,1]$ acting as
integration operator with kernel $\sigma^2$.

Moreover, Gaussian random matrices with a certain degree of correlation between the
entries are also asymptotically semicircular elements over an appropriate subalgebra, see
\cite{ROBS}.

\sect{Further Free-Probabilistic Aspects of Random Matrices}

Free probability theory provides also new ideas and techniques for investigating other
aspects of random multi-matrix models. In particular, Haagerup and Thorbjornsen
\cite{Haa02,Haa05} obtained a generalization to several matrices for a number of results
concerning the largest eigenvalue of a Gaussian random matrix.

Much work is also devoted to deriving rigorous results about the large $N$ limit of
random multi-matrix models given by densities of the type
$$c_N e^{-N^2 \tr P(A_1,\dots,A_N)}d\lambda(A_1,\dots,A_n),$$
where $d\lambda$ is Lebesgue measure, $A_1,\dots,A_n$ are selfadjoint $N\times N$
matrices, and $P$ a noncommutative selfadjoint polynomial. To prove the existence of that
limit in sufficient generality is one of the big problems. For a mathematical rigorous
treatment of such questions, see \cite{Gui06}.

\emph{Free Brownian motion} is the large $N$ limit of the Dyson Brownian motion model of
random matrices (with independent Brownian motions as entries, compare Chapter 11). Free
Brownian motion can be realized concretely in terms of creation and annihilation
operators on a full Fock space (see section \ref{sect:freegroup}). There exists also a
corresponding \emph{free stochastic calculus} \cite{Bia98}; for applications of this to
multi-matrix models, see \cite{Gui07}.

There is also a surprising connection with the representation theory of the symmetric
groups $S_n$. For large $n$, representations of $S_n$ are given by large matrices which
behave in some respects like random matrices. This was made precise by Biane who showed
that many operations on representations of the symmetric group can asymptotically be
described by operations from free probability theory, see \cite{Bia02}.

\sect{Operator Algebraic Aspects of Free Probability}

A survey on free probability without the mentioning of at least some of its operator
algebraic aspects would be quite unbalanced and misleading. We will highlight some of
these operator algebraic facets in this last section. For the sake of brevity, we will
omit the definitions of standard concepts from operator algebras, since these can be
found elsewhere (we refer the reader to \cite{VDN,Voi05,Hia00} for more information on
notions, as well as for references related to the following topics).

\subsection{Operator Algebraic Models for Freeness}\label{sect:freegroup}

\subsubsection{Free group factors} Let $G=\star_{i\in I} G_i$ be the free product of groups
$G_i$. Let $L(G)$ denote the group von Neumann algebra of $G$, and $\ff$ the associated
trace state, corresponding to the neutral element of the group. Then $L(G_i)$ can be
identified with a subalgebra of $L(G)$ and, with respect to $\ff$, these subalgebras
$(L(G_i))_{i\in I}$ are free. This freeness is nothing but the rewriting in terms of
$\ff$ what it means that the groups $G_i$ are free as subgroups in $G$. The definition of
freeness was modeled according to the situation occurring in this example. The
\emph{free} in free probability theory refers to this fact.

A special and most prominent case of these von Neumann algebras are the \emph{free group
factors} $L(\FF_n)$, where $\FF_n$ is the free group on $n$ generators. One hopes to eventually
be able to resolve the isomorphism problem: whether the free groups factors $L(\FF_n)$ and
$L(\FF_m)$ are, for $n,m\geq 2$, isomorphic or not.

\subsubsection{Creation and annihilation operators on full Fock spaces}

Let $\HH$ be a Hilbert space. The \emph{full Fock space} over $\HH$ is defined as
$\cF(\HH):=\bigoplus_{n=0}^\infty \HH^{\otimes n}$. The summand $\HH^{\otimes 0}$ on the
right-hand side of the last equation is a one-dimensional Hilbert space. It is customary to
write it in the form $\CC\Omega$ for a distinguished vector of norm one, which is called
\emph{the vacuum vector}. The vector state $\tau_{\HH}$ on $B(\cF(\HH))$ given by the vacuum
vector, $\tau_\HH(T):=\langle T \Omega, \Omega\rangle$ ($T\in B( \cF(\HH))$), is called
\emph{vacuum expectation state}.

For each $\xi \in \HH$, the operator $l( \xi ) \in B( \cF ( \HH ))$ determined by the
formula $l( \xi ) \Omega = \xi$ and $l(\xi)\xi_1\otimes\dots\otimes \xi_n = \xi \otimes
\xi_1\otimes\dots\otimes \xi_n$ for all  $n \geq 1$, $\xi_1,\ldots , \xi_n \in \HH$, is
called the \emph{(left) creation operator} given by the vector $\xi$. As one can easily
verify, the adjoint of $l( \xi )$ is described by the formula: $l( \xi )^* \Omega = 0$,
$l( \xi )^* \xi_1 = \langle \xi_1,\xi \rangle \Omega$, and $l( \xi )^*
\xi_1\otimes\dots\otimes \xi_n =
 \langle \xi_1 ,\xi \rangle \xi_2 \otimes\dots\otimes \xi_n$
and is called the \emph{(left) annihilation operator} given by the vector $\xi$.

The relevance of these operators comes from the fact that orthogonality of vectors translates
into free independence of the corresponding creation and annihilation operators.

\begin{proposition}
Let $\HH$ be a Hilbert space and consider the probability space $(B(\cF(\HH)),\tau_\HH)$. Let
$\HH_1, \ldots , \HH_k$ be a family of linear subspaces of $\HH$, such that $\HH_i \perp \HH_j$
for $i \neq j$ $(1 \leq i,j \leq k).$ For every $1 \leq i \leq k$ let $\cA_i$ be the unital
$C^*$-subalgebra of $B( \cF ( \HH ))$ generated by $\{ l ( \xi ) : \xi \in \HH_i \}$. Then
$\cA_1,\dots,\cA_k$ are freely independent in $(B(\cF(\HH)),\tau_\HH)$.
\end{proposition}

Also semicircular elements show up very canonically in this frame; namely, if we put
$l:=l(\xi)$ for a unit vector $\xi\in\HH$, then $l+l^*$ is a semicircular element of variance
1. More generally, one has that $l+f(l^*)$ has $\cR$-transform $\cR(z)=f(z)$ (for $f$ a
polynomial, say). This, together with the above proposition, was the basis of Voiculescu's
proof of Theorem \ref{thm:R-transform}. Similarly, a canonical realization for the
$S$-transform is $(1+l)g(l^*)$, for which one has $S(z)=1/g(z)$. This representation is due to
Haagerup who used it for a proof of Theorem \ref{thm:S-transform}.

\subsection{Free Entropy}

Free entropy is, as the name suggests, the counterpart of entropy in free probability
theory. The development of this concept is at present far from complete. The current
state of affairs is that there are two distinct approaches to free entropy. These should
give isomorphic theories, but at present we only know that they coincide in a limited
number of situations. The first approach to a theory of free entropy is via
\emph{microstates}. This goes back to the statistical mechanics roots of entropy via the
Boltzmann formula and is related to the theory of large deviations. The second approach
is \emph{microstates free}. This draws its inspiration from the statistical approach to
classical entropy via the notion of Fisher information. We will in the following only
consider the first approach via microstates, as this relates directly with random matrix
questions.

Wigner's semicircle law states that as $N \rightarrow \infty$ the empirical eigenvalue
distribution $\mu_{A_N}$ of an $N\times N$ Gaussian random matrix $A_N$ converges almost
surely to the semicircular distribution $\mu_W$, i.e., the probability that $\mu_{A_N}$
is in any fixed neighborhood of the semicircle converges to 1. We are now interested in
the deviations from this: What is the rate of decay of the probability that $\mu_{A_N}$
is close to $\nu$, where $\nu$ is an arbitrary probability measure? We expect that this
probability behaves like $e^{-N^2I(\nu)}$, for some \emph{rate function} $I$ vanishing at
the semicircle distribution. By analogy with the classical theory of large deviations,
$I$ should correspond to a suitable notion of free entropy. This heuristics led
Voiculescu to define in \cite{Voi93} the \emph{free entropy} $\chi$ in the case of one
variable to be
$$\chi(\nu)=\iint \log\vert s-t\vert d\nu(s)d\nu(t)+\frac 34+\frac 12 \log 2\pi.$$
Inspired by this, Ben-Arous and Guionnet proved in \cite{Ben97} a rigorous version of a
large deviation for Wigner's semicircle law, where the rate function $I(\nu)$ is, up to a
constant, given by $-\chi(\nu)+ \frac 12 \int t^2 d\nu(t)$.

Consider now the case of several matrices. By Voiculescu's generalization of Wigner's
theorem we know that $n$ independent Gaussian random matrices $A^{(1)}_N,\dots,
A^{(n)}_N$ converge almost surely to a freely independent family $s_1,\dots,s_n$ of
semicircular elements. Similarly as for the case of one matrix, large deviations from
this limit should be given by
\begin{equation*}
\text{Prob}\left\{(A^{(1)}_N,\dots,A^{(n)}_N) : \distr((A^{(1)}_N,\dots,A^{(n)}_N)
\approx \distr(a_1,\dots,a_n)\right\} \sim e^{-N^2I(a_1,\dots,a_n)},
\end{equation*}
where $I(a_1,\dots,a_n)$ should be related to the \emph{free entropy} of the random
variables $a_1,\dots,a_n.$ Since the distribution $\distr(a_1,\dots,a_n)$ of several
non-commuting random variables $a_1,\dots,a_n$ is a mostly combinatorial object
(consisting of the collection of all joint moments of these variables), it is much harder
to deal with these questions and, in particular, to get an analytic formula for $I$.
Essentially, the above heuristics led Voiculescu to the following definition
\cite{Voi94b} of a free entropy for several variables.

\begin{definition}
Given a tracial $W^*$-probability space $(M,\tau)$ (i.e., $M$ a von Neumann algebra and
$\tau$ a faithful and normal trace), and an $n$-tuple $(a_1,\dots,a_n)$ of selfadjoint
elements in $M$, put
\begin{align*}
\Gamma(a_1,\dots,&a_n;N,r,\epsilon):= \\ &\bigl\{(A_1,\dots,A_n) \in M_N(\CC)_{sa}^n:
|\tr(A_{i_1}\dots A_{i_k})-\tau(a_{i_1}\dots a_{i_k})|\leq \epsilon\\
&\qquad\qquad\qquad\qquad\qquad\qquad\text{ for all } 1 \leq i_1,\dots,i_k \leq n, 1 \leq
k \leq r\bigr\}
\end{align*}
In words, $\Gamma(a_1,\dots,a_n;N,r,\epsilon)$ is the set of all $n$-tuples of $N \times
N$ selfadjoint matrices which approximate the mixed moments of the selfadjoint elements
$a_1,\dots,a_n$ of length at most $r$ to within $\epsilon.$

Let $\Lambda$ denote Lebesgue measure on $M_N(\CC)_{sa}^n.$  Define
\begin{equation*}
\chi(a_1,\dots,a_n;r,\epsilon):=\limsup_{N\to\infty} \frac{1}{N^2} \log
\Lambda(\Gamma(a_1,\dots,a_n;N,r,\epsilon)+ \frac{n}{2}\log N,
\end{equation*}
and
\begin{equation*}
\chi(a_1,\dots,a_n):=\lim_{\substack{r \rightarrow \infty\\
\epsilon \rightarrow 0}} \chi(a_1,\dots,a_n;r,\epsilon).
\end{equation*}
The function $\chi$ is called the \emph{free entropy}.
\end{definition}

Many of the expected properties of this quantity $\chi$ have been established (in
particular, it behaves additive with respect to free independence), and there have been
striking applications to the solution of some old operator algebra problems. A celebrated
application of free entropy was Voiculescu's proof of the fact that free group factors do
not have Cartan subalgebras (thus settling a longstanding open question). This was
followed by several results of the same nature; in particular, Ge showed that $L({\mathbb
F}_n)$ cannot be written as a tensor product of two $\text{II}_1$ factors. The rough idea
of proving the absence of some property for the von Neumann algebra $L({\mathbb F}_n)$
using free entropy is the following string of arguments: finite matrices approximating in
distribution any set of generators of $L({\mathbb F}_n)$ should also show an approximate
version of the considered property; one then has to show that there are not many finite
matrices with this approximate property; but for $L({\mathbb F}_n)$ one has many
matrices, given by independent Gaussian random matrices, which approximate its canonical
generators.

However, many important problems pertaining to free entropy remain open. In particular,
we only have partial results concerning the relation to large deviations for several
Gaussian random matrices. For more information on those and other aspects of free entropy
we refer to \cite{Voi02, Bia03, Gui04}.

\subsection[Operator Algebraic Applications]{Other Operator Algebraic Applications of Free Probability
Theory} The fact that freeness occurs for von Neumann algebras as well as for random
matrices means that the former can be modeled asymptotically by the latter and this
insight resulted in the first progress on the free group factors since Murray and von
Neumann. In particular, Voiculescu showed that a compression of some $L(\FF_n)$ results
in another free group factor; more precisely, one has
$(L(\FF_n))_{1/m}=L(\FF_{1+m^2(n-1)})$. By introducing interpolated free group factors
$L(\FF_t)$ for all real $t>1$, this formula could be extended by Dykema and Radulescu to
any real $n,m>1$, resulting in the following dichotomy: One has that either all free
group factors $L(\FF_n)$ $n\geq 2$ are isomorphic or that they are pairwise not
isomorphic.

There exist also type III versions of the free group factors; these free analogues of the
Araki-Woods factors were introduced and largely classified by Shlyakhtenko.

The study of free group factors via free probability techniques has also had an important
application to subfactor theory. Not every set of data for a subfactor inclusion can be
realized in the hyperfinite factor; however, work of Shlyakhtenko, Ueda, and Popa has
shown that this is possible using free group factors.

By relying on free probability techniques and ideas, Haagerup achieved also a crucial
break-through on the famous invariant subspace problem: every operator in a $II_1$ factor whose
Brown measure (which is a generalization of the spectral measure composed with the trace to
non-normal operators) is not concentrated in one point has non-trivial closed invariant
subspaces affiliated with the factor.

\ \\
{\sc Acknowledgements}: This work was supported by a Discovery Grant from NSERC.


\begin{thebibliography}{Abc84a}

\bibitem[And10]{And10}
G. Anderson, A. Guionnet, and O. Zeitouni, {\it An Introduction to Random Matrices}, Cambridge
University Press (to appear)

\bibitem[Ben97]{Ben97}
G. Ben-Arous and A. Guionnet, Prob. Th. Rel. Fields {\bf 108} (1997) 517

\bibitem[Ben09]{Ben09}
F. Benaych-Georges,  Prob. Th. Rel. Fields {\bf 144} (2009) 471

\bibitem[Ber93]{Ber93}
H. Bercovici and D. Voiculescu, Indiana Univ. Math. J. {\bf 42} (1993) 733

\bibitem[Ber99]{Ber99}
H. Bercovici and V. Pata (with an appendix by P. Biane), Ann. of Math. {\bf 149} (1999)
1023

\bibitem[Bia98]{Bia98}
P. Biane and R. Speicher, Prob. Th. Relat. Fields {\bf 112} (1998) 373

\bibitem[Bia02]{Bia02}
P. Biane, Proceedings of the International Congress of Mathematicians, Beijing 2002, Vol.
2 (2002) 765

\bibitem[Bia03]{Bia03}
P. Biane, M. Capitaine, and A. Guionnet, Invent. Math. {\bf 152} (2003) 433

\bibitem[Cvi82]{Cvi82} P. Cvitanovic, P.G.~Lauwers, and P.N.~Scharbach, Nucl. Phys. B {\bf
203} (1982) 385

\bibitem[Col07]{Col07}
B. Collins, J. Mingo, P. Sniady, and R., Speicher, Documenta Math. {\bf 12} (2007) 1

\bibitem[Dyk93]{Dyk93}
K. Dykema, J. Funct. Anal. {\bf 112} (1993) 31

\bibitem[Gui04]{Gui04} A.~Guionnet, Probab. Surv. {\bf 1} (2004) 72

\bibitem[Gui06]{Gui06}
A. Guionnet, Proceedings of the International Congress of Mathematicians, Madrid 2006,
Vol. III (2006) 623

\bibitem[Gui07]{Gui07}
A.~Guionnet and D.~Shlyakhtenko, preprint, arXiv:math/0701787 (2007)

\bibitem[Haa02]{Haa02}
U. Haagerup, Proceedings of the International Congress of Mathematicians, Beijing 2002,
Vol 1 (2002) 273

\bibitem[Haa05]{Haa05}
U. Haagerup, S. Thorbj\o rnsen, Ann. of Math. {\bf 162} (2005) 711

\bibitem[Hia00]{Hia00} F. Hiai and D. Petz, {\it The Semicircle Law,
Free Random Variables and Entropy}, Math. Surveys and Monogr. 77, AMS 2000

\bibitem[Nic06]{NS}
A. Nica and R. Speicher, {\it Lectures on the Combinatorics of Free Probability}, London
Mathematical Society Lecture Note Series, vol. 335, Cambridge University Press, 2006

\bibitem[Min10]{MS} J. Mingo and R. Speicher, {\it Free Probability and Random Matrices},
Fields Monograph Series (to appear)

\bibitem[Rao09]{Rao09}
N.~R.~Rao and A. Edelman, Foundations of Computational Mathematics (to appear)

\bibitem[Ras08]{ROBS}
R.~Rashidi Far, T.~Oraby, W.~Bryc, and R. Speicher, IEEE Trans. Inf. Theory {\bf 54}
(2008) 544

\bibitem[Shl96]{Shl} D. Shlyakhtenko, IMRN 1996 {\bf 1996} 1013

\bibitem[Spe93]{Spe93}
R. Speicher, Publ. RIMS {\bf 29} (1993) 731

\bibitem[Spe94]{Spe94}
R. Speicher, Math. Ann. {\bf 298} (1994) 611

\bibitem[Spe98]{Spe98}
R. Speicher, {\it Combinatorial theory fo the free product with amalgamtion and
operator-valued free probability theory}, Memoirs of the AMS {\bf 627} 1998

\bibitem[Tul04]{Tul04} A. Tulino, S. Verdu, {\it Random matrix
theory and wirless communications}, Foundations and Trends in Communications and
Information Theory {\bf 1} (2004)

\bibitem[Voi85]{Voi85}
D.~Voiculescu, in {\it Operator Algebras and their Connections with Topology and Ergodic
Theory}, Lecture Notes in Math. 1132 (1985), Springer Verlag, 556

\bibitem[Voi91]{Voi91}
D. Voiculescu, Invent. Math. {\bf 104} (1991) 201

\bibitem[Voi92]{VDN}
D. Voiculescu, K. Dykema, and A. Nica, {\it Free Random Variables}, CRM Monograph Series,
Vol. 1, AMS 1992

\bibitem[Voi93]{Voi93}
D. Voiculecu, Comm. Math. Phys. {\bf 155} (1993) 71

\bibitem[Voi94a]{Voi94}
D. Voiculescu, Proceedings of the International Congress of Mathematicians, Z\"urich
1994, 227

\bibitem[Voi94b]{Voi94b}
D. Voiculescu, Invent. Math. {\bf 118} (1994) 411

\bibitem[Voi95]{Voi95}
D. Voiculescu, Asterisque {\bf 223} (1995) 243

\bibitem[Voi00]{Voi00}
D. Voiculescu, in {\it Lectures on Probabiltiy Theory and Statistics (Saint-Flour,
1998)}, Lecture Notes in Mathematics 1738, Springer, 2000, 279

\bibitem[Voi02]{Voi02}
D. Voiculescu, Bulletin of the London Mathematical Society {\bf 34} (2002) 257

\bibitem[Voi05]{Voi05}
D. Voiculescu, Reports on Mathematical Physics {\bf 55} (2005) 127

\bibitem[Woe86]{Woe86}
W. Woess, Bollettino Un. Mat. Ital. {\bf 5-B} (1986) 961

\bibitem[Xu97]{Xu97}
F. Xu, Commun. Math. Phys. {\bf 190} (1997) 287

\end{thebibliography}
\end{document}